\documentclass[12pt,reqno]{amsart} 
\usepackage{amsthm,amsfonts,amssymb,amscd}

\newtheorem{theorem}[subsection]{Theorem}
\newtheorem{proposition}[subsection]{Proposition}

\newtheorem{lemma}[subsection]{Lemma}
\newtheorem{corollary}[subsection]{Corollary}

\theoremstyle{definition}

\newtheorem{proposition-definition}[subsection]{Proposition-Definition}

\theoremstyle{remark}
\newtheorem{remark}[subsection]{Remark}

\newcommand{\Hth}{\stackrel{\circ}{H}_{3,0}}
\newcommand{\Hff}{H_{4,0}}
\newcommand{\Hfo}{\stackrel{\circ}{H}_{4,0}}

\newcommand{\corank}{\operatorname{corank}\nolimits}

\newcommand{\Hom}{{{\mathcal H}om\:}}
\newcommand{\Homg}{\operatorname{Hom}\nolimits}
\newcommand{\id}{\operatorname{id}\nolimits}
\newcommand{\Ext}{\operatorname{Ext}\nolimits}
\newcommand{\Alb}{\operatorname{Alb}\nolimits}

\newcommand{\Tr}{\operatorname{Tr}}
\newcommand{\elm}{\operatorname{elm}\nolimits}
\newcommand{\Gr}{\operatorname{Gr}}
\newcommand{\Pic}{\operatorname{Pic}\nolimits}
\newcommand{\pr}{\operatorname{pr}\nolimits}

\newcommand{\End}{{{\mathcal E}nd\:}}

\newcommand{\PGL}{{PGL}_5(\CC )}

\newcommand{\CC}{{\mathbb C}}

\newcommand{\ZZ}{{\mathbb Z}}

\newcommand{\PP}{{\mathbb P}}

\newcommand{\FF}{{\mathbb F}}

\newcommand{\OOO}{{\mathcal O}}
\newcommand{\OOP}{{\mathcal O}_{\PP^1}}
\newcommand{\III}{{\mathcal I}}

\newcommand{\EEE}{{\mathcal E}}
\newcommand{\HHH}{{\mathcal H}}

\newcommand{\FFF}{{\mathcal F}}
\newcommand{\CCC}{{\mathcal C}}

\newcommand{\NNN}{{\mathcal N}}
\newcommand{\MMM}{{\mathcal M}}
\newcommand{\EXT}{{\mathcal Ext}}

\newcommand{\EEEE}{\boldsymbol{\mathcal E}}
\newcommand{\GGGG}{\boldsymbol{\mathcal G}}

\newcommand\epsi{\epsilon}

\newcommand\g{\gamma}

\newcommand\lra{{\longrightarrow}}
\newcommand\rar{\rightarrow}
\newcommand\Hilb{{\operatorname{Hilb}\nolimits^{5n}_X}}
\newcommand\HILB[1]{{\operatorname{Hilb}\nolimits^{#1}_X}}
\newcommand\Hilbp{{\operatorname{Hilb}\nolimits^{5n}_{\PP^4}}}
\newcommand\HILBP[1]{{\operatorname{Hilb}\nolimits^{#1}_{\PP^4}}}
\newcommand\Sym{{\operatorname{Sym}\nolimits}}
\newcommand\ns{{\operatorname{ns}\nolimits}}

\renewcommand\square{\frame{\phantom{{\large x}}}}

\author{D. Markushevich}

\address{D. M.: Math\'ematiques - b\^{a}t. M2, Universit\'e Lille 1,
F-59655 Villeneuve d'Ascq Cedex, France}
\email{markushe@gat.univ-lille1.fr}

\author{A. S. Tikhomirov}

\thanks{The work supported by the Volkswagen-Stiftung (RiP
program at Oberwolfach)}

\address{A. T.: Department of Mathematics, State Pedagogical University,
Respublikanskaya 108, Yaroslavl' 150000, Russia}
\email{alexandr@tikho.yaroslavl.su}


\title{The Abel--Jacobi map of a moduli component of vector bundles
on the cubic threefold}

\begin{document}
\begin{abstract}
The Abel--Jacobi map of the family of elliptic quintics lying
on a general cubic threefold is studied. It is proved that
it factors through a moduli component of stable rank 2 vector
bundles on the cubic threefold with Chern numbers $c_1=0,
c_2=2$,
whose general point represents
a vector bundle obtained by Serre's construction from an elliptic
quintic. The elliptic quintics mapped to a point of the moduli
space vary in a 5-dimensional projective space
inside the Hilbert scheme of curves,
and the map from the moduli space to the intermediate Jacobian
is \'etale.  As auxiliary results, the irreducibility of families
of elliptic normal quintics and of rational normal quartics
on a general cubic threefold is proved. This implies the uniqueness
of the moduli component under consideration. The techniques of
Clemens--Griffiths and Welters are used for 
the calculation of the infinitesimal Abel--Jacobi map.
\end{abstract}
\maketitle

\section*{Introduction}

In their famous paper \cite{CG}, Clemens and Griffiths represented
the intermediate Jacobian $J_1(X)$ of a smooth cubic
threefold $X$ as the Albanese variety of the family
of lines on $X$ and parametrized
the theta divisor on it by the Abel--Jacobi images of pairs of 
lines. They used this parametrization to establish
the Torelli Theorem for $X$ and the non-rationality of $X$. 
Similar results were
obtained by Tyurin \cite{Tyu}. Later several authors constructed
parametrizations of the intermediate Jacobians and, sometimes, 
of their theta divisors for
other Fano threefolds, using families of curves of low degree lying on them,
see, for example, \cite{CV}, \cite{De-1}, \cite{De-2},
\cite{I}, \cite{Lo}, \cite{Le},
\cite{T-1}, \cite{T-2}, \cite{V-1}, \cite{V-2}, \cite{We}.

However, the following natural question has not been investigated:
given a family of connected curves parametrized by
some variety $B$, what are the fibers of the Abel--Jacobi map
considered on $B$ itself
and not on the $\Alb (B)$? By analogy with the Abel--Jacobi map
from divisors of a given degree $d$ on a curve of 
genus $g$ to its Jacobian, which is smooth with projective
spaces as fibers, provided that $d\geq 2g-1$, 
we might expect that the Abel--Jacobi map in our context behaves
well for families of curves of rather big degree. The
papers mentioned above search, on the contrary, for 
the families of the smallest posssible degree whose image generates the
intermediate Jacobian. For the cubic threefold $X$, the only known
facts in this direction concern the map of difference, defined on
the pairs of lines in $X$, which is
onto the theta divisor and is generically of degree 6 (see \ref{psi}), and the
Abel--Jacobi map of the family of rational cubics in $X$, which is
also onto the theta divisor and its generic fiber is $\PP^2$ (see Proposition
\ref{prop41}). There are no similar results about the Abel--Jacobi
parametrizations of the entire intermediate Jacobian.

In the present paper, we study the Abel--Jacobi
map on the 10-dimensional component $H$
of the Hilbert scheme whose general point represents an elliptic
normal quintic in $X$. We find an open subset $\HHH\subset H$, on which
the map is smooth and all the components of its fibers are isomorphic
to $\PP^5$ (Theorem \ref{main}). The analogy with 
the Jacobians of curves is not
complete, because the fibers may consist of {\em many}
copies of the projective
space, and almost nothing is known about the boundary of $\HHH$ in $H$
(see Proposition \ref{aftermain}).

We give also an unexpected interpretation to the variety $M=M_X$
parametrizing the connected components of the fibers of the
Abel--Jacobi map $\Phi : \HHH\lra J_1(X)$: we show that it is
isomorphic to an open subset in the component of the moduli space
$M_X(2;0,2)$ of stable rank 2 vector
bundles on $X$ with Chern numbers $c_1=0,
c_2=2$, whose general member is obtained by Serre's construction
from elliptic quintics in $X$. Thus, there exists a natural
Abel--Jacobi map from this component of the moduli space 
of stable vector bundles on $X$ to $J_1(X)$, induced by $\Phi$;
this map is nothing but the Abel--Jacobi map of the moduli
space, defined by the second Chern class of the vector bundle
with values in the Chow group of 1-cycles modulo rational equivalence.
Moreover,
this map is \'etale on the open subset $M$, see Theorem \ref{main}.

As far as we know, the above theorem provides the first example
of a moduli space of vector bundles which has a dominant map to an
abelian variety, different from the Picard and Albanese varieties
of the base. It also shows that the moduli spaces of vector bundles
yield sometimes more efficient parametrizations of the intermediate
Jacobian of a threefold than families of curves.

This construction can be relativized over the family of nonsingular
hyperplane sections of a cubic fourfold $V$. It gives a
moduli component of torsion sheaves on $V$
with supports on the hypeplane sections $X=H\cap V$, whose restrictions
to $X$ are rank two vector bundles from $M_X$. One can show that
the Yoneda pairing induces a symplectic structure on it. This is
a new example of a moduli space of sheaves possessing a symplectic
structure; the known ones parametrize sheaves on hyperkaehler 
manifolds (see \cite{Muk}
for surfaces, and \cite{K} for higher dimensional varieties).
It is a 10-dimensional
symplectic variety, covering the one constructed by Donagi
and Markman in \cite[Example 8.22]{D-M} (work in process).

Now, we will briefly describe the contents of the paper by sections.

In Section 1, we remind basic facts about the cubic threefolds,
the Abel--Jacobi map and its differential, following essentially
\cite{CG} and \cite{We}.

Section 2 explains Serre's construction, applied to elliptic
normal quintics in $X$. It yields a 5-dimensional component of
the moduli space of stable vector bundles of rank 2 on $X$.
The unstability of vector bundles obtained by Serre's construction
from {\em non-normal} elliptic quintics is also proved. 

Section 3 gives the proof of the irreducibility
of the family of rational normal quartics $\Gamma\subset X$. First, 
the irreducubility of the family of rational twisted cubics 
lying in nonsingular hyperplane sections of $X$ is
proved, in using the Abel--Jacobi map of these curves to the
theta divisor in $J_1(X)$. Next, the irreducibility of the family
of curves of the form $D+l$ is proved, where $D$ is a twisted cubic as above
and $l$ a line meeting $D$ transversely at 1 point. We show also that
the Hilbert scheme $\HILB{4n+1}$ is smooth at points representing
such curves $D+l$, and that they are strongly smoothable into
a rational normal quartic (the techniques
of \cite{HH} are used). This implies the wanted result.

Section 4 is devoted to the normal elliptic
quintics $C\subset\PP^4$ contained in a cubic 
threefold $X$. It is proved that the
Hilbert scheme $\Hilb$ is smooth at points representing such curves $C$,
and that their family is irreducible
for a {\em general} $X$. The proof is reduced to that of
the irreducibility of the family of rational normal quartics $\Gamma\subset X$,
in using the liaison defined by cubic scrolls $\FF_1\subset \PP^4$:
$\FF_1\cap X=C\cup\Gamma$.

Section 5 contains the main result of the paper. We relativize
Serre's construction of Section 2 over some open
subset $H_0\subset H$, and construct a morphism $\phi :H_0\lra M_X(2;0,2)$
with fibers $\PP^5$. We prove that $\phi$ is, locally in the
\'etale topology, the structure projection of a projectivized
vector bundle of rank 6. The smoothness of the Abel--Jacobi map
$\Phi :\HHH\lra J_1(X)$ is proved in using the technique of Tangent
Bundle Sequence, and this, together with the fact that
an Abel--Jacobi map should contract to a point any projective space,
implies that the components of fibers of $\Phi$ are exactly the fibers
of $\phi$, isomorphic to $\PP^5$. We also study the extension $\tilde{\Phi}$
of $\Phi$ to the boundary component of $\HHH\subset H$, formed by non-normal
elliptic quintics, and find that its image coincides with a translate
of $F+F\subset \Alb (F)$, where $F$ is the Fano surface of $X$,
and that the differential of $\tilde{\Phi}$ is degenerate along this
component.

\smallskip

{\em Acknowledgements}. We acknowledge with pleasure the hospitality
of the Mathematisches Forschungsinstitut at Oberwolfach, where
we made an essential part of the work during our RiP stay. 
We are grateful to J.~Harris
for his idea of the link between
ellipitic quintics and rational quartics on a cubic threefold, and
to A.~Beauville for his remark, cited in \ref{phi}.
We are also grateful to K.~Hulek, A.~Iliev, E.~Markman,
S.~Mukai, C. Peskine, and F.-O.~Schreyer for discussions.

\section{Generalities and known results}

\subsection{}
We will start by a reminder on cubic threefolds. Let $X$ be a smooth
cubic hypersurface in $\PP^4$. It is a Fano variety, that is its
anticanonical sheaf $\omega_X^{-1}\simeq\OOO_X(2)$ is ample,
and the following properties are obtained by standard techniques for Fano threefolds
\cite{Isk}:
\begin{equation}\label{generalities-1}\begin{array}{c}
h^i(\OOO_X(k))=0\; \mbox{for}\; i=1,2, k\in\ZZ\; , \\ 
\; h^{i,0}=h^{0,i}=0
\; \mbox{for}\; i>0, h^{1,2}=h^{2,1}=5\; ,
\end{array}
\end{equation} 
\begin{equation}\label{generalities-2}\begin{array}{c}
\Pic (X)=A_2(X)=H^2(X,\ZZ )=\ZZ \cdot [\OOO_X(1)]\; , \\
B_1(X)=H^4(X,\ZZ )=\ZZ\cdot l \; .\end{array}
\end{equation}
Following \cite[Ch. 19]{Fu}, we denote by $A_i$, resp. $B_i$ the Chow group of 
algebraic cycles of dimension $i$ modulo rational, resp.
algebraic equivalence, and $l$ is the class of a line. 

\subsection{}
The geometry of lines (and, sometimes, of conics) plays a crucial role 
in the study of a Fano threefold $X$, essentially by the following
two reasons. Firstly, they give (a part of) generators of the group
of birational automorphisms of $X$ \cite{Isk-2}, and secondly, they are
used for the description of the Abel--Jacobi map
$AJ:\Homg_1X\lra J_1(X)$, where we denote, following again \cite{Fu},
by $\Homg_iX$ the subgroup in the group 
$Z_iX$ of algebraic $i$-cycles, consisting of cycles 
homologous to 0, and $J_1(X)$ is the intermediate
Jacobian of $X$. According to \cite{CG}, the lines on a nonsingular
cubic threefold $X$ are parametrized by a nonsingular surface
$F=F(X)\subset \Gr (2,5)$ with invariants
\begin{equation}\label{hijF}
h^{1,0}(F)=5\; ,\; h^{2,0}(F)=10,
\end{equation}
which can be thought of as a component of the Hilbert scheme
$\operatorname{Hilb}\nolimits^{n+1}_X$. The smoothness of
$\operatorname{Hilb}\nolimits^{n+1}_X$ at the points of $F$
follows from the calculation of the normal bundle of a line
$l\subset X$:
\begin{equation}\label{NlX}
\begin{array}{c}
\NNN_{l/X}\simeq \OOO\oplus\OOO\;\; (\mbox{for}\; l\in F\setminus D) \\
\NNN_{l/X}\simeq \OOO (-1)\oplus\OOO (1)\;\; (\mbox{for}\; l\in
D), \end{array}\end{equation}
where $D\subset F$ is a curve. Formulas (\ref{NlX}) imply that
$h^1(\NNN_{l/X})=0$, hence $F$ is smooth by \cite{G}. $F$ is called the
{\em Fano surface} of $X$. 

\subsection{}\label{AJ_B}
We will recall some known facts about
the Abel--Jacobi map \cite{Gri}, \cite{CG}, \cite{Tyu}, \cite{We}.
The intermediate Jacobian of a threefold $X$ is defined by
$$
J_1(X)=(F^2H^3(X,\CC ))^*/im\: (H_3(X,\ZZ )),
$$
where $F^2=H^{3,0}+H^{2,1}$ is a term of the Hodge filtration, and
$H_3(X,\ZZ )$ is mapped to $(F^2H^3(X,\CC ))^*$ by integration over
cycles. For a cubic threefold, $H^{3,0}(X)=0$, so, by the Hodge
index theorem, the intersection form on $H_3(X)$ defines a
principal polarization on $J_1(X)$. 

The Abel--Jacobi map is defined as follows: let $\gamma\in\Homg_1X$.
Then there exists $\Gamma\in H_3(X,\ZZ )$ such that $\gamma =
\partial\Gamma$, and $AJ(\gamma )$ is given by
$$
\left[ \omega\mapsto\int_\Gamma \omega\right] \mod im\: H_3(X,\ZZ ).
$$
For a family of 1-cycles $\{ Z_b\}$ parametrized by some variety  $B$ with a
reference point $\beta$,
this defines a set-theoretic map $AJ_B:B\lra J_1(X), \;
b\mapsto AJ(Z_b-Z_\beta )$. By \cite[II]{Gri}, $AJ_B$ is analytic
if $B$ is nonsingular. Hence $AJ_B$ factors
through the Albanese variety of (the desingularization of) $B$;
the resulting morphism $\Alb (B)\lra J_1(X)$ is also
called Abel--Jacobi map.
Applying this construction to the family
of lines on a cubic threefold $X$, Clemens--Griffiths get a natural
map $\Alb (F)\lra J_1(X)$.  By (\ref{generalities-1})
and (\ref{hijF}), both abelian varieties are of dimension 5.
The following facts are known:

\subsection{}
The Abel--Jacobi map establishes an isomorphism of 
abelian varieties $\Alb (F)\simeq J_1(X)$ (\cite[(0.9)]{CG}).

\subsection{}
For any fixed $s\in F$, the map $i_s:F\lra J_1(X)$, $t\mapsto [l_t-l_s]$ is 
a closed embedding, where $l_t$ denotes the line in $X$ corresponding to
a point $t\in F$ (\cite[\S 3]{Tyu}).

\subsection{}\label{theta3}
The class of $i_s(F)$ in $B_2(J_2(X))$ is $\frac{1}{3!}\theta^3$,
where $\theta$ is the class of the theta divisor defining the principal
polarization on $J_2(X)$ (\cite[(0.9)]{CG}).

\subsection{}\label{psi}
The map $\psi :F\times F\lra J_1(X)$, $(s,t)\mapsto [l_t-l_s]$ is
generically 6-to-1 over the image, and $[\psi (F\times F)]=\theta$
(ibid, (0.10)). 

A precise description of the singularities of the map $\psi$ is given
in \cite{Tyu}.
See also loc. cit., or \cite[Example 4.3.2]{Fu} and references therein
for further information on the geometry of the Abel--Jacobi map of
lines.

\subsection{}\label{phi}\label{beau}
We will use in the sequel another map $\tilde{\psi} :F\times F\lra J_1(X)$,
$(s,t)\mapsto [l_s+l_t-2l_0]$, where $0\in F$ is some reference point.
The same arguments as in \cite{CG} show that $\tilde{\psi} $ is generically
finite over its image and \ref{theta3} implies that its class in
$B_2(J_1(X))$ is $m\theta$ for some $m\in\ZZ, m>0$.  
The generical finiteness follows from the tangent bundle theorem
for $F$ \cite[12.4]{CG}, which implies that whenever $l_s,l_t$
is a pair of skew lines, the tangent spaces to $i_0(F)$ in
$\Alb (F)$ at $i_0(s),i_0(t)$ are transversal (when translated to
$0\in\Alb (F)$), so that both maps of sum ($\tilde{\psi} $) and of difference
($\psi$) have injective differentials at $(s,t)$.
Beauville communicated us that the values of
the degree $d$ of the map $\tilde{\psi} :F\times F\lra \tilde{\psi} (F\times F)$
and of $m$ follow from the description of the intermediate Jacobian
in terms of Prym varieties \cite{B}: $d=2$ and $m=3$.

\subsection{Technique of TBS}\label{TBS}
Now we will describe, following \cite[Sect. 2]{We}, the technique
of ``tangent bundle sequence", which will be applied later to
elliptic quintics in a cubic threefold.

Let $X\hookrightarrow W$ be an embedding of a smooth projective threefold
$X$ into a smooth quasiprojective fourfold $W$, and $\{ Z_b\}_{b\in B}$
a flat family of curves in $X$ parametrized by a nonsingular
variety $B$ over $\CC$. Then a choice of a base point $\beta\in
B$ determines the Abel--Jacobi map $B\lra J_1(X)$. Let $Z=Z_b$ be
a scheme theoretical locally complete intersection fiber of the family.
Then the differential of the Abel--Jacobi map at $b$ factors into 
the composition of the two natural maps:
$$
T_{B,b}\lra H^0(Z,\NNN_{Z/X})
$$
and
$$
\psi_Z:H^0(Z,\NNN_{Z/X})\lra (H^0(X,\Omega^3_X)\oplus
H^1(X,\Omega^2_X))^*=T_{J_1(X),0}
$$
The $H^0(X,\Omega^3_X)^*$-component of $\psi_Z$ is always 0, so in the
sequel we will consider and denote by the same symbol $\psi_Z$ the map to
$H^1(X,\Omega^2_X)^*$. In \cite{We}, this factorization is given for
$Z$ smooth with a reference to \cite[II]{Gri}. This holds also in a more
general context \cite{Fl}: if we identify $H^0(\NNN_{Z/X})$ with
$\Ext^1_X(\FFF ,\FFF)$, where $\FFF =\OOO_Z$, then $\psi_Z$ can be 
represented as the composition $\Ext^1_X(\FFF ,\FFF)
\raisebox{0.5 ex}{$\begin{CD}
@>{\cup at(\FFF )}>> \end{CD}$}
\Ext^2_X(\FFF ,\FFF\otimes\Omega^1 )
\stackrel{\Tr}{\lra} H^2(\Omega^1)$,
where $at(\FFF )\in\Ext^1_X(\FFF ,\FFF\otimes\Omega^1 )$ denotes the
Atiyah--Illusie class of a coherent sheaf $\FFF$.

It is convenient to describe $\psi_Z$ via its conjugate $\psi_Z^*$,
which fits into the following commutative square:
\begin{equation}\label{CDWelters}
\begin{CD}
H^0(X,\NNN_{X/W}\otimes\omega_X) @>{R}>> H^1(X,\Omega^2_X) \\
@V{r_Z}VV @VV{\psi_Z^*}V \\
H^0(Z,\NNN_{X/W}\otimes\omega_X|_Z) @>{\beta_Z}>> H^0(Z, \NNN_{Z/X})^*.\\
\end{CD}
\end{equation}
Here $r_Z$ is the map of restriction to $Z$, and the whole square
(upon natural identifications)
is a part of the commutative diagram of long exact cohomology
sequences associated to the following commutative diagram
of sheaves:
\begin{equation}\label{CDsheaves}
\begin{array}{ccccccccc}
\scriptstyle{ 0} &\scriptstyle{ \rar} &\scriptstyle{ \Omega^2_X\otimes\OOO_Z} &\scriptstyle{ \rar}
&        \scriptstyle{ \Omega^3_W\otimes\NNN_{X/W}\otimes\OOO_Z }&\scriptstyle{ \rar}
&                  \scriptstyle{  \Omega^3_X\otimes\NNN_{X/W}\otimes\OOO_Z }&\scriptstyle{ \rar }&\scriptstyle{ 0 }\\
& & \downarrow & &\downarrow & & \| & & \\
\scriptstyle{ 0} &\scriptstyle{\rar} &\scriptstyle{ \Omega^3_X\otimes \NNN_{Z/X} }&\scriptstyle{ \rar} &\scriptstyle{ \Omega^3_X\otimes\NNN_{Z/W}}
      &  \scriptstyle{    \rar} &\scriptstyle{ \Omega^3_X\otimes\NNN_{X/W}\otimes\OOO_Z} &\scriptstyle{ \rar} &\scriptstyle{ 0}
\end{array}
\end{equation}
See \cite[2.8]{We} for definitions of the maps.

\section{Vector bundles from elliptic quintics}

\label{sect.5}

In the sequel, $X$ will denote a nonsingular cubic
threefold in $\PP^4$,
$H$ or $H(X)$ the union of components of $\Hilb$
having smooth normal elliptic quintics $C\in\PP^4$ as their generic points.
We will see later in Section \ref{section-quintics} that $H(X)$ is non empty
and has dimension 10 for any $X$.
We will denote by $[C]_V$ the point representing a subscheme $C\subset V$
in the Hilbert scheme of $V$; the subscript $V$ may be omitted, if
this makes no confusion.

In this section, we will study the vector bundles $\mathcal E$ on $X$ obtained
by Serre's construction applied to an elliptic quintic $C$ in $X$.
They are defined by the following extension of $\OOO_X$-modules:
\begin{equation}\label{serre}
0\lra \OOO_X\lra \EEE (1) \lra \III_C(2) \lra 0\; ,
\end{equation}
where $\III_C=\III_{C,X}$ is the ideal sheaf of $C$ in $X$.
Since the class of $C$ in $B_1(X)$
is $5l$, the sequence (\ref{serre}) implies that 
$c_1(\EEE )=0, c_2(\EEE )=2l$. Moreover, $\det\EEE$ is trivial,
and hence  $\EEE$ is self-dual as soon as it is a vector
bundle (that is, $\EEE^\ast\simeq\EEE$).

We will verify that there exists a unique non-trivial
extension (\ref{serre}) for given $C$ (Corollary \ref{unique-ext}).

\begin{lemma}\label{dim=6} For any smooth $C\in H$, we have: 

a) $\dim\Ext^1(\III_C(2),\OOO_X )=1$;

b) for $k\geq 2,\; h^0(\III_C(k))=\binom{4+k}{4}-\binom{1+k}{4}
-5k,\; h^1(\III_C(k))=$
$h^2(\III_C(k))=0$; $h^0(\III_C(1))=h^1(\III_C(1))$ is $0$ if $C$ is not
contained in a hyperplane, and $1$ otherwise, and $h^2(\III_C(1))=0$;

c) $h^0(\III_C (2))=5,h^0(\EEE (1))=6,h^i(\EEE (1))=0$ for $i>0$.
\end{lemma}

\begin{proof}
a) Applying $\Ext (\cdot ,\OOO_X)$ to the restriction sequence
\begin{equation}\label{restriction}
0\lra\III_C(k)\lra\OOO_X(k)\lra\OOO_C(k)\lra 0 
\end{equation}
for $k=2$, we obtain:
\begin{equation}\label{Ext}\begin{array}{r}
\Ext^1(\OOO_X (2),\OOO_X)\lra\Ext^1(\III_C(2),\OOO_X)\lra
\Ext^2(\OOO_C (2),\OOO_X)\\
\lra\Ext^2(\OOO_X (2),\OOO_X) 
\end{array}
\end{equation}
By (\ref{generalities-1}), the left and the right hand terms
are zero. So, $\Ext^1(\III_C(2),\OOO_X)=
\Ext^2(\OOO_C (2),\OOO_X)=H^0(X,\EXT_{\OOO_X}^2(\OOO_C,\omega_X))=
H^0(C,\omega_C)\simeq\CC$.

b) Looking again at (\ref{restriction}), we see that the assertion
for $k\geq 1$ follows from the projective normality of $C$
\cite[Proposition IV.1.2]{Hu}
in the case when $C$ does not lie in a hyperplane.

Let now $C\subset\PP^3$. If $k=1$,
the assertions follow easily from (\ref{restriction}).
For $k= 2$,
it is sufficient to show that $C$ is not contained in a quadric.
Assume the contrary. $C$ is obviously not contained in a nonsingular
quadric, neither in the nonsingular locus of a quadratic cone,
because in that case its degree should be even. So it lies in a quadratic
cone $Q\subset\PP^3$ and passes through the vertex. By degree
considerations, it has to meet twice the generators of the cone 
outside the vertex. So, in blowing up the vertex,
we obtain a nonsingular curve $\tilde{C}$ of genus 1 in the Hirzebruch
surface $\FF_2$ with class $2s+af$, where $s,f$ are the standard generators
of $\Pic \FF_2$ such that $s^2=-2, f^2=0, (s\cdot f)=1$. Computing
the canonical class of $\tilde{C}$, we find immediately $a=4$, hence
$(s\cdot \tilde{C})=0$, which implies that $C$ does not pass through
the vertex. This contradicts our assumptions. Hence $C$ does not lie
in any quadric.

We have seen that $h^1(\III_C (2))=h^2(\III_C (1))=0$; obviously, also
$h^3(\III_C)=0$, so that $C$ is 3-regular according to the definition
of Castelnuovo--Mumford. The Castenuovo's Proposition in Lecture 14
of \cite{Mum} implies the case $k>2$.

c) The first equality follows from b) with $k=2$.
The remaining ones are immediate consequences of (\ref{serre})
and of the Serre duality.

\end{proof}

\begin{corollary}\label{unique-ext}
For any smooth  $C\in H$, 
there exists a unique (up to isomorphism)
extension (\ref{serre}) with $\EEE$ locally free.
\end{corollary}

\begin{proof}
The uniqueness follows from Lemma \ref{dim=6}, a). The local freeness
follows from the sheafified version of (\ref{Ext}), giving
$\EXT_{\OOO_X}^1(\III_C(2),\OOO_X)=
\EXT_{\OOO_X}^2(\OOO_C,\omega_X)=\omega_C$,
and from the following lemma due to Serre:

\begin{lemma}[{\cite[Ch. I, Lemma 5.1.2]{Okonek}}] 
Let $X$ be a nonsingular variety, 
$C$ a locally complete intersection of codimension 2 in $X$, 
$L$ invertible on $X$, and
$$
0\lra \OOO_X\lra \EEE  \lra \III_C \otimes L \lra 0\
$$
an extension given by a class  $e\in \Ext^1(\III_C\otimes L,\OOO_X)$.
Then $\EEE$ is locally free of rank 2 if and only if the image of $e$ in
$H^0(C,\EXT_{\OOO_X}^1(\III_C\otimes L,\OOO_X))$ generates the stalk
of $\EXT_{\OOO_X}^1(\III_C\otimes L,\OOO_X)$ at every point of $C$.
\end{lemma}
\end{proof}

\begin{corollary}\label{P5toH}
If $C$ is smooth and is not contained in a hyperplane, then
any non-zero section of $\EEE (1)$ has a locally complete intersection curve
with zero canonical class
as its zero locus, and this defines a natural morphism
$\PP^5=\PP (H^0(X,\EEE (1))\lra H$ whose image contains $[C]$.
\end{corollary}

\begin{proof}
The scheme of zeros of a section of a rank 2 vector bundle 
on a nonsingular variety is locally
complete intersection as soon as it is of codimension 2;
the assertion about the canonical class follows by adjunction.
So, we have to show that a non-zero section of $\EEE (1)$
cannot vanish on a surface $S\subset X$. By (\ref{generalities-2}),
$\OOO_X(S)\simeq\OOO_X(d)$ for some $d>0$, hence it suffices to show that
$h^0(\EEE (1-d))= 0$ for all $d> 0$. Twisting (\ref{serre}) by
$\OOO_X(-1)$, we see that if $C$ is not contained in a hyperplane,
that is, $h^0(\III_C(1))=0$, then $h^0(\EEE )=0$, and so much the more
$h^0(\EEE (1-d))= 0$. This ends the proof.
\end{proof}

For future use in Section \ref{factor}, we will study also the case
of non-normal elliptic quintics.

\begin{proposition}\label{space-C}
If $C$ is a smooth elliptic quintic contained in a hyperplane section
$S$ of $X$, then there exists a unique pair of possibly coincident lines
$l_1\cup l_2$ in $S$, such that the zero locus of any section of
$\EEE (1)$, if not of codimension 2, is of the form 
$S'\cup l_1\cup l_2$ for some hyperplane section $S'$ of $X$.
In this case $C$ is rationally equivalent to a curve of the
form $C^3+l_1+l_2$, where $C^3$ is a plane cubic in $S$.

Conversely, any curve of the form $C^3+l_1+l_2$ on a general cubic
threefold $X$, where $C^3$ is a plane cubic curve and $l_1,l_2$
a pair of disjoint lines, is rationally equivalent to a
smooth elliptic quintic contained in a hyperplane section of $X$.
\end{proposition}
\begin{proof}
As in the previous corollary, we find a section $s$ of $\EEE (1)$
vanishing on $S$, and a section $s_0$ of $\EEE$ such that $s=Ls_0$,
where $L$ is the linear form defining $S$ in $X$. As $c_2( \EEE )=2l$,
the scheme of zeros $Z$ of $s_0$ is either a conic, or a pair of
lines; $Z$ may be non-reduced (a double line), but without embedded
points, for it is a locally complete intersection. As $h^0(\EEE (1))=6$,
we deduce from the exact triple
\begin{equation}\label{nonstable}
0\lra \OOO_X(1)\stackrel{\alpha}{\lra} \EEE (1)\stackrel{\beta}{\lra} \III_Z(1)\lra 0
\end{equation}
that $h^0(\III_Z(1))=1$, so that $Z$ is contained in a unique
hyperplane. This eliminates the case of the conic, as well as that
of a planar double line.
Thus, $\chi(\OOO_Z(n))=2n+2,\ <Z>=P^3$, i.~e. either $Z$ is a disjoint union 
$Z=l_1\bigsqcup l_2$ of two lines, or $Z$ is a double non-planar locally complete
intersection structure on a line $l_1=l_2$ (the latter is clearly the limit 
of the former).
We will denote this shortly as $Z=l_1\cup l_2$.
Moreover, we see from (\ref{nonstable}) that
$\EEE (1)$ has a 5-dimensional subspace of sections
$L's_0$ with $L'\in H^0(X,\OOO_X(1))$ whose zero locus is
of the form $S'\cup l_1\cup l_2$, where $S'=\{ L'=0\}$ is a
hyperplane section of $X$. 

Let $s_1$ be the section of $\EEE (1)$ vanishing exactly on $C$,
and $L'$ a general linear form. Then the pencil of sections
$\lambda_0L's_0+\lambda_1s_1$ defines a pencil of curves, degenerating
$C$ into $C^3+l_1+l_2$ with $C^3\subset S'$ of degree 3. 

Now, $s_1$ is not in the image of $\alpha$ in (\ref{nonstable}),
and it vanishes on $C$. Hence $\beta (s_1)\neq 0$ and $\beta (s_1)$
vanishes on $C$. Hence $C$ lies in the unique hyperplane
containing $Z=l_1\cup l_2$. Hence, this hyperplane is $L=0$, and
$l_1\cup l_2\subset S$.

It remains to see that $C^3\subset S$. By construction, $C^3\subset S'$
for another hyperplane section $S'$ of $X$. If $C^3$ were not contained in
$S\cap S'$, then the general member of the pencil of curves defined
by the sections $\lambda_0L's_0+\lambda_1s_1$ would be a smooth elliptic
curve not contained in a hyperplane, because $C^3\cup l_1\cup l_2$
is not. This contradicts Corollary \ref{P5toH}. Hence $C^3
\subset S\cap S'$ is a plane cubic curve.

The fact that any curve of the form $C^3+l_1+l_2$ with $l_1\cap
l_2 =\varnothing$ is rationally equivalent
to a smooth quintic is reduced to the study of linear systems on
cubic surfaces. We can always arrange the things so that $C^3,l_1,l_2$
lie in one hyperplane and $C^3$ meets each one of the lines $l_1,l_2$
transversely at 1 point, because all the plane sections $C^3$ of $X$
are rationally equivalent to each other, being parametrized by
the  rational Grassmannian variety $G(3,5)$. We obtain
in this way the linear system $|C^3+l_1+l_2|$ on the cubic surface
$S=<l_1,l_2>\cap X$. If $S$ is nonsingular, we can represent it as $\PP^2$ 
with 6 points $P_1,\ldots ,P_6$ blown  up. Let $e_0$ be the inverse image
of a line in $\PP^2$, and $e_i$ the exceptional divisors over $P_i$
($i=1,\ldots ,6$). We can choose this representation in such a way,
that $l_1=e_1,l_2=e_2$. We have $C^3\in |h|$, where $h=3e_0-e_1-\ldots -e_6$ 
is the hyperplane section, so that $C^3+l_1+l_2\in |3e_0-e_3-\ldots -e_6|$
is obviously smoothable in its linear system.

To treat the case when $S$ is singular, we will use the assumption that
$X$ is general. This restricts the possible degenerations of the
hypeplane sections $H\cap X$ to those which have codimension $\leq 4$
in the projective space $\PP^{19}$ parametrizing cubic surfaces in
$\PP^3$. Such degenerations $S$ have at most 4 isolated singular points,
and the line joining two distinct singular points lies entirely in $S$
(see \cite[Sect. 4.6]{GH}).
Choosing coordinates in such a way that the singular points lie in
the vertices of the coordinate octahedron and writing out the monomials that
can be present in the equation of $S$ in order that $S$ might have
isolated singularities in the prescribed vertices,
we arrive at the following list of degenerations up to codimension 4:
$A_1$, $A_1A_1$, $A_2$, $A_1A_1A_1$, $A_2A_1$, $A_3$, $A_1A_1A_1A_1$, $
A_2A_1A_1$, $A_3A_1$, $A_2A_2$, $D_4$,
where $X_i$ (resp. $X_iX_j,\ldots$) stands for a surface with
one singular point of type $X_i$ (resp. distinct singular points
of types $X_i,X_j,\ldots$), and the analytic types of singularities $X_i$ are
either $A_i$ or $D_i$ of Du Val, defined up to analytic equivalence
by $x^2+y^2+z^{i+1}=0\; (A_i)$ or $x^2+y^2z+z^{i-1}=0\; (D_i)$.
The codimension of such degenerations in $\PP^{19}$ is 
$i$ (resp. $i+j+\ldots$).

The standard technique of projections from two lines (see
\cite{M} or \cite{Re}) gives a birational map $S\dasharrow
\PP^1\times\PP^1$. In the above 11 cases, it represents $S$ as
$\PP^1\times\PP^1$ with 5 points $\Psi=\{ Q_1,\ldots, Q_5\}$
blown up and ($-2$)-curves blown down; some of the points $Q_i$
can be infinitesimal. Since the quadric with a point blown up
is isomorphic to $\PP^2$ with two points blown up, we can replace
$(\PP^1\times\PP^1,\Psi )$ by $(\PP^2,\Phi )$, where
$\Phi =\{ P_1, \ldots , P_6\}$ is a set of 6 points, some of
which can be infinitesimal. Let $S_\Phi$ denote the blowup of
$\Phi$ in $\PP^2$.

We are going to provide a configuration of $\Phi$ giving rise to
each of the 11 degenerations above. In each case, $\Phi$ should
be general satisfying the formulated constraints in the sense
that there are no curves $F\simeq\PP^1, F^2\leq -2$ among the
proper transforms in $S_\Phi$ of lines or conics in $\PP^2$
except for the ($-2$)-curves imposed by the constraints. Here
is the list:
\begin{itemize}

\item[$A_1$:] $P_1,P_2,P_3$ are collinear.

\item[$A_1A_1$:] $P_3=P_1P_2\cap P_4P_5$.

\item[$A_2$:] $P_1,P_2,P_3$ are collinear, and
$P_4,P_5,P_6$ are collinear.

\item[$A_1A_1A_1$:] $P_2,P_4,P_6$ lie on the sides of the triangle
$P_1P_3P_5$.

\item[$A_2A_1$:] The limit of $A_1A_1$ when $P_4\to P_2$ along the line
$P_2P_5$ (so that $P_4$ is infinitesimally close to $P_2$).

\item[$A_3$:] $P_3=P_1P_2\cap P_4P_5$, and $P_6$ is infinitesimally close
to $P_3$).

\item[$A_1A_1A_1A_1$:] $\Phi$ is the set of intersection points of 4 lines
in general position in~$\PP^2$.

\item[$A_2A_1A_1$:] The specialization of $A_2A_1$ with $P_6\in P_1P_5$.

\item[$A_3A_1$:] The limit of $A_2A_1$ when $P_6\to P_1$.

\item[$A_2A_2$:] The limit of $A_2A_1$ when $P_6\to P_5$.

\item[$D_4$:] $P_1,P_2,P_3$ are collinear, and $P_{3+i}$ is
infinitesimally close to $P_i\; (i=1,2,3)$.
\end{itemize}

A routine verification of the emptiness of the base  loci in $S_\Phi$
of the linear systems $h+l_1+l_2$ ends the proof.

\end{proof}

\begin{proposition}\label{stable}
Let $C$ be a smooth elliptic quintic. If $C$ is not 
contained in a hyperplane, then
the vector bundle $\EEE$ obtained from $C$ by Serre's construction
is stable. It has a non-trivial global section and is Gieseker unstable, 
if $C$ is contained in a hyperplane.
\end{proposition}
\begin{proof}
{\em Unstability part}. Assume that $C$ is contained in a hyperplane.
As in the proof of Proposition \ref{space-C},
we obtain an exact triple (\ref{nonstable}). It is obvious that
$\chi (\OOO_X(k))>\chi (\III_Z(k))$ for $k\gg 0$, so $\EEE$
is Gieseker unstable, though semistable in the sense of Mumford--Takemoto.

{\em Stability part}. Any saturated torsion free rank 1 subsheaf of $\EEE(1)$ 
is invertible of the form $\OOO_X(k)$ and gives an exact triple:
\begin{equation}\label{diese}
0\lra \OOO_X(k)\stackrel{\alpha}{\lra} \EEE (1)\lra \III_{Z}(2-k)
\lra 0
\end{equation}
where $Z\subset X$ is a subscheme of codimension 2 (if not empty).
Clearly, (\ref{serre}) twisted by $\OOO_X(-2)$ shows that the case $k\ge2$ here
is impossible. On the other hand, if $\EEE$ is not Gieseker stable,
then, computing Hilbert polynomials, one has $k\ge1$ in (\ref{diese}).
Hence $k=1$, and we obtain the above case when $C$ lies in a hyperplane.
\end{proof}

\begin{lemma}\label{h0NCX}
Let $C$ be a smooth elliptic quintic in $X$. Let $\EEE$ be the
unique locally free sheaf defined by (\ref{serre}). Then we have:

a) $h^i(\EEE (-1))=0\;\; \forall\; i\in\ZZ$;

b) $h^0(\EEE\otimes\EEE )=1$ if and only if $C$ is not contained in a hyperplane,
and in this case $h^1(\EEE\otimes\EEE )=5,
h^2(\EEE\otimes\EEE )=h^3(\EEE\otimes\EEE )=0$;

c) $h^0(\NNN_{C/X})=10,h^1(\NNN_{C/X})=0$.
\end{lemma}

\begin{proof}
a) From (\ref{restriction}) with $k=0$, we deduce $h^i(\III_C)=0,\ i=0,1,3,$
$h^2(\III_C)=1$. Then (\ref{serre}), twisted by
$\OOO_X(-2)$, gives $\chi(\EEE (-1))=$ \linebreak 
$h^i(\EEE (-1))=0$ for $i=0,1$. 
Next, by Serre duality,  $h^3(\EEE (-1))=h^0(\EEE (-1))=0$; hence
also $h^2(\EEE (-1))=0$.

b) Tensor (\ref{serre}) by $\EEE (-1)$:
\begin{equation}\label{EE}
0\lra \EEE (-1)\lra \EEE\otimes\EEE\lra \III_C\otimes\EEE (1)\lra 0
\end{equation}
>From a) and (\ref{EE}), we deduce that 
\begin{equation}\label{hi}
h^i(\EEE\otimes\EEE)=h^i(\III_C\otimes\EEE (1))
\;\;\forall\; i\in\ZZ .
\end{equation}
If $C$ is not contained in a hyperplane, then
$\EEE$ is stable, hence by \cite[Theorem I.2.9]{Okonek}, $\EEE$
is simple, that is  $h^0(\Hom (\EEE ,\EEE ))=h^0(\EEE\otimes\EEE )=1$.
By Serre duality, $h^3(\EEE\otimes\EEE)=h^0(\EEE\otimes\EEE (-2))=0$.
By Riemann--Roch--Hirzebruch, $\chi (\EEE\otimes\EEE)=-4$.
It remains to prove that $h^2(\EEE\otimes\EEE)=h^1(\EEE\otimes\EEE
(-2))=0$.

Consider the exact triples
$$
0\lra \EEE (-3)\lra \EEE\otimes\EEE (-2)\lra \III_C\otimes\EEE (-1)\lra 0
$$
$$
0\lra \III_C\otimes\EEE (-1)\lra \EEE (-1)\lra \EEE (-1)\otimes\OOO_C\lra 0
$$
We have $\EEE (-1)\otimes\OOO_C\simeq\NNN_{C/X}(-2)$.
Consider the natural inclusion 
$\NNN_{C/X}(-2)\subset\NNN_{C/\PP^4}(-2)$ and apply
\cite[Proposition V.2.1]{Hu}: \linebreak 
$h^1(\NNN^\ast_{C/\PP^4}(2)\otimes \MMM)=0$
for any invertible $\MMM$ of degree 0. By Serre duality, 
$h^0(\NNN_{C/\PP^4}(-2))=0$.
Hence $h^0(\NNN_{C/X}(-2))=0$. The last exact sequence and
a) imply that $h^i(\III_C\otimes\EEE (-1))=
h^i(\EEE (-1))=0\;,\; i=0,1$, hence
$h^1(\EEE\otimes\EEE (-2))=h^1(\EEE (-3))$.
By Lemma \ref{dim=6}, c), and Serre duality, we have $h^2(\EEE (1))
=h^1(\EEE (-3))=0$, and we are done.

If $C$ is contained in a hyperplane, then $h^0(\EEE )\neq 0$ by
Proposition \ref{stable}, hence $\EEE\otimes\EEE=\End (\EEE )$ has
at least two linearly independent sections: the identity endomorphism,
and $s\otimes s$, where $s$ is a non-trivial section of $\EEE$.

c) Let $C$ be not contained in a hyperplane.
In the exact triple
\begin{equation}\label{timesE1}
 0\lra \III_C\otimes\EEE (1)\lra \EEE (1)\lra \NNN_{C/X}\lra 0
\end{equation}
we have $h^0( \III_C\otimes\EEE (1))=1, h^1( \III_C\otimes\EEE (1))=5$
by b) and (\ref{hi}), and $h^0(\EEE(1))=6,h^1(\EEE(1))=0$ by
Lemma \ref{dim=6}, c). Hence $h^0(\NNN_{C/X})=10$. By Riemann--Roch
on $C$, we have  $h^1(\NNN_{C/X})=0$.

In the case if $C\subset\PP^3$, we can consider $C$ as a curve
on a surface, the hyperplane section $S$ of $X$, as in the
proof of Proposition \ref{space-C}.
For example, assume that the hyperplane section
is nonsingular. Then we can represent it as $\PP^2$ with 6 points 
$P_1,\ldots ,P_6$ blown  up. Let $e_0$ be the inverse image of a line in
$\PP^2$, and $e_i$ the exceptional divisors over $P_i$ ($i=1,\ldots ,6$).
Then $C\in |3e_0-e_{i_1}-e_{i_2}-e_{i_3}-e_{i_4}|$ for some subset
$\{ i_1,i_2,i_3,i_4\}\subset\{1,\ldots 6\}$.
The standard exact sequences
$$
0\lra \OOO_S\lra \OOO_S(C)\lra \NNN_{C/S}\lra 0,
$$
$$
0\lra \NNN_{C/S}\lra \NNN_{C/X}\lra \OOO_C(1)\lra  0
$$
show that $h^0(\NNN_{C/X})=10,h^1(\NNN_{C/X})=0$.

\end{proof}

\begin{lemma} Let $C$ be a smooth elliptic quintic,  not 
contained in a hyperplane.
Then any two non-proportional  sections of $\EEE (1)$ define different
curves $C$. Hence the morphism $\PP^5\lra H$ of Corollary \ref{P5toH}
is injective.
\end{lemma}

\begin{proof}
By Lemma \ref{h0NCX}, b),  $h^0(\EEE\otimes\EEE )=1$. To prove
that the section of $\EEE (1)$ defining $C$ is unique 
up to a constant factor,
it suffices  to show that $h^0( \III_C\otimes\EEE (1))=1$.
This follows from (\ref{hi}). 
\end{proof}

\section{Rational normal quartics}

In this section, $X$ is a general cubic threefold. For the future use
in Section \ref{section-quintics},
we prove the irreducibility of the family of rational normal quartics
in $X$.

\begin{theorem}\label{quart-irred}
There is a unique component of the Hilbert 
scheme of curves in $X$ whose generic
point correponds to a rational normal quartic.
\end{theorem}

We start by several auxiliary results.

\begin{proposition}\label{prop41}
Let $X$ be a nonsingular cubic threefold, and $\Hth =$ 
$\Hth (X) $ the open subscheme of $\HILB{3n+1}$
parametrizing smooth cubic rational curves lying in nonsingular
hyperplane sections of $X$. Then we have:

(i) $\Hth$ is smooth and irreducible of dimension 6.

(ii) The image of the Abel--Jacobi map
$\eta :\Hth\lra J_1(X)$ is an open subset of a translate $\Theta_a$ of
the theta divisor of $J_1(X)$.

(iii) All the scheme-theoretic fibers $\eta^{-1}\eta (s)$ are open
subschemes of $\PP^2$, in particular, they are smooth and irreducible
of dimension 2.
\end{proposition}

\begin{proof}
See \cite[(6.21)]{We}. We have
\begin{equation}\label{Csmooth}
\eta^{-1}\eta (s) = \{ C\in |\OOO_S(C_s)|\; : \; C \;\mbox{is smooth}
\},
\end{equation}
where $C_s\subset X$ is the curve represented by $s\in\Hth$,
$S=<C_s>\cap X$, and $<C_s>\simeq\PP^3$ denotes the linear span
of $C_s$. Thus, the $\PP^2$ from the statement of the proposition
is the linear system of $C_s$ on the nonsingular cubic surface $S$.
The image of $\eta$ can be identified, up to translation, with
the image under $\psi$ (see \ref{psi}) of the subset
$U=\{ (t,t')\in F\times F\; |\; l_t\cap l_{t'}=\varnothing ,\
<l_t, l_{t'}>\cap \ X\ \mbox{is nonsingular}\}$.
The pairs $(t,t')$ are recovered from $s$ as follows:
there are exactly 6 lines $l_t$ in $S$ such that $C\in |\OOO_S(C_s)|$ degenerates 
into $C^2+l_t$, where $C^2$ is a conic. Choose one of them, and let
$l_{t'}$ be determined as the residual intersection
of the linear span of $C^2$ with $X$: $X\cap <C^2>=C^2\cup l_{t'}$.
Then the Abel--Jacobi image of $C_s$ coincides, up to a translation
by a constant, with that of $l_t-l_{t'}$.
The openness of $\eta (\Hth )$ follows from the finiteness of $\psi$
on $U$ \cite[\S 3]{Tyu}.
\end{proof}

\begin{proposition}\label{prop42}
Let $X$ be a general cubic threefold, $\Hth$ as in Proposition
\ref{prop41}, and $F$ the Fano surface of $X$. Let 
$$
B=B(X)=\{ (s,t)\in\Hth\times F\; |\; C_s \;\mbox{\em meets}\; l_t\;
\mbox{\em transversely at 1 point}\} ,
$$
where $C_s$ (resp. $l_t$) denotes the rational cubic curve
represented by $s$ (resp. the line represented by $t$). Then
$B$ is irreducible of dimension $7$.
\end{proposition}

\begin{proof}
Let $b:B\lra\Theta_a$ be defined by $b=\eta\circ p_1$, where $\Theta_a,
\eta$ are as in Proposition \ref{prop41}, and $p_1$ is the natural projection
from $B$ to $\Hth$. The map $b$ is dominant on every component
of $B$. To see this, notice first that on a {\em general} cubic 
$X$, the number of lines passing through any point is finite ($\leq 6$).
Indeed, by \cite[\S 1]{Tyu} or \cite[Sect. 8]{CG}, the points $P\in X$
lying on an infinite number of lines are characterized by the property that
$T_PX\cap X$ is a cone over a plane cubic curve, and it is easy to check by
counting dimensions that there are no such hyperplane sections on
a general $X$.

Hence {\em all} the fibers $p_1^{-1}(x)$ are of pure dimension 1. Moreover,
any $(x_0,t_0)\in p_1^{-1}(x_0)$ can be obtained as a limit of points
$(x,t)\in p_1^{-1}(x)$ for some one-parameter family specializing $x$
into $x_0$. Hence $B$ is equidimensional and each one of its components
dominates $\Hth$. This implies, by Proposition \ref{prop41}, that each
component of $B$ dominates $\Theta_a$.

Look now at the fibers of $b:B\lra\Theta_a$. Let $v\in\eta (\Hth )
\subset\Theta_a$, $B_v=b^{-1}(v)$, and $f=p_2|_{B_v}:B_v\lra F$,
$f:([C],[l])\mapsto [l]$.
By (\ref{Csmooth}),
$$\begin{array}{r}
B_v=\{ ([C'],[l])\in\Hth\times F\; |\; C'\in|\OOO_S(C)|,
\; l\;\mbox{meets}\; C'\;\mbox{transversely}\phantom{oint}\\ 
\mbox{in 1 point}\} .\end{array}
$$
Hence, if $l\not\subset S$,
$$
f^{-1}([l])=\{\mbox{open subset of}\;\PP^1=|\III_{x,S}(C)|\} ,
$$
where $x=l\cap S$. So, all the fibers of $f:B_v\lra F$ are open
subsets of $\PP^1$ except for 27 fibers $\simeq\PP^2$, corresponding to
lines $l_i\subset S$ ($i=1,\ldots ,27$). These exceptional fibers are not
irreducible components, as any $l_i\cup C'$ represented by a point of
$B_v$ can be deformed into $l\cup C''$ with $l\not\subset S$,
$C''\cap l=\{1 \;\mbox{point}\}$, in using the fact that the family of
rational cubics covers $S$. Hence $B_v$ is irreducible of constant dimension
3, hence $B$ is irreducible, and we are done.
\end{proof}

\begin{proposition}\label{prop43}
Let $X,B$ be as in Proposition \ref{prop42}. Then $B$ can be
naturally identified with a subset of $\Hff =\HILB{4n+1}$.
The following assertions are true:

(i) $\Hff$ is smooth at any point $[Z]=([C],[l])\in B$, and
$$
\dim T_{[Z]}\Hff =h^0(\NNN_{Z/X})=8,\; h^1(\NNN_{Z/X})=0 .
$$

(ii) A curve $Z=C\cup l$ represented by a general point of $B$ can be
smoothed into a rational normal quartic inside $X$.
\end{proposition}

\begin{proof}
(i) We have natural exact sequences
\begin{equation}\label{seq1}
0\lra\NNN_{Z/X}\lra \NNN_{Z/X}|_C\oplus \NNN_{Z/X}|_l
\lra \NNN_{Z/X}\otimes\CC (x)\lra 0,
\end{equation}
\begin{equation}\label{seq2}
0\lra \NNN_{C/X}\lra \NNN_{Z/X}|_C\lra \CC (x)\lra 0,
\end{equation}
\begin{equation}\label{seq3}
0\lra \NNN_{l/X}\lra \NNN_{Z/X}|_l\lra \CC (x)\lra 0,
\end{equation}
where $x=C\cap l$.

Let $S=<C>\cap X$, then $\NNN_{C/S}\simeq\OOP (1)$ via an identification
$C\simeq\PP^1$. Hence the exact sequence
$$
0\lra \NNN_{C/S}\lra\NNN_{C/X}\lra \NNN_{S/X}|_C\lra 0
$$
yields
\begin{equation}\label{normcubic}
\NNN_{C/X}\simeq \OOP (2)\oplus\OOP (2)\;\mbox{or}
\; \OOP (1)\oplus\OOP (3).
\end{equation}

By (\ref{NlX}),(\ref{normcubic}), (\ref{seq2}) and  (\ref{seq3}),
$h^1(\NNN_{C/X})=h^1(\NNN_{l/X})=0$,
$h^0(\NNN_{Z/X}|_C)=7$, $h^0(\NNN_{Z/X}|_l)=3$. The map
$H^0(\NNN_{Z/X}|_C\oplus \NNN_{Z/X}|_l)
\lra H^0(\NNN_{Z/X}\otimes\CC (x))$ arising from (\ref{seq1})
is surjective, because $\NNN_{Z/X}$ is locally free
and $\NNN_{Z/X}|_C\simeq \OOP (2)\oplus\OOP (3)\;\mbox{or}
\; \OOP (1)\oplus\OOP (4)$ is generated by global sections.
Hence $h^0(\NNN_{Z/X})=8$,
$h^1(\NNN_{Z/X})=0$.

(ii) We will use the following particular case of
\cite[Theorem 4.1]{HH} (though it is formulated in \cite{HH} for nodal
curves in $\PP^3$, its proof remains valid with $\PP^3$ replaced by
any nonsingular projective variety).

\begin{lemma}
Let $X$ be a nonsingular projective variety, $C_1,C_2$ two smooth
curves in $X$ meeting transversely at one point $x$. Assume that
$H^1(C_i,\elm^-_x\NNN_{C_i/X})=H^1(C_{3-i},\NNN_{C_{3-i}/X})=0$
for at least one $i=1$ or $2$. Then $H^1(\NNN_{C_1\cup C_2/X})=0$
and $C_1\cup C_2$ is strongly smoothable in $X$.
\end{lemma}
Here $\elm^-_x$ denotes the negative elementary transformation
in $x$ (see loc. cit., Sect. 2 for a definition).

>From (\ref{normcubic}), we obtain
$$
\elm^-_x\NNN_{C/X}\simeq\OOP (1)\oplus\OOP (2)\;\mbox{or}
\; \OOP \oplus\OOP (3).
$$
Hence $h^1(\elm^-_x\NNN_{C/X})=0$. We have already seen that
$h^1(\NNN_{l/X})=0$, hence the lemma implies the result.
\end{proof}

\subsubsection*{Proof of Theorem \ref{quart-irred}}
Let $\Hff '(X)$ be the component of $\Hff (X)$ containing $B(X)$,
where $B(X)$ is defined in Proposition \ref{prop42}.
It is unique, \hfill because \hfill $\Hff (X)$ 
 is smooth and of dimension 8 at the points of
 $B(X)$,  and thus  we  have  $\dim \Hff '(X)=8$.

Let $\Hfo\subset\HILBP{4n+1}$ be the subscheme parametrizing
rational normal quartics; it is a smooth homogeneous manifold of dimension
21: $\Hfo\simeq PGL_5(\CC )/PGL_2(\CC )$. Let
$$
I=\{ (X,C)\in |\OOO_{P^4}(3)|\times\Hfo\; |\; C\subset X\}
$$
be the incidence variety, $p_1,p_2$ its projections to the factors
of \linebreak $|\OOO_{P^4}(3)|\times\Hfo$. Then
$p_2^{-1}([C])=|\III_{C,P^4}(3)|\simeq\PP^{21}$ for any $C$.
Hence $I$ is irreducible of dimension 42. We have seen that there is
an irreducible component of a general fiber of $p_1:I\lra\PP^{34}
=|\OOO_{P^4}(3)|$ of dimension 8, hence $p_1$ is dominant. Moreover,
the irreducibility of $I$ implies that the action of the monodromy
group on the ireducible components of a general fiber is transitive.
As we can distinguish one particular component $\Hff '(X)$, invariant
under monodromy, the general fiber of $p_1$ is irreducible.
\hfill\square

\begin{remark}
The basic idea of our proof of the irreducibility of the family
of rational normal quartics in $X$ is degenerating a quartic into
a reducible curve, which reduces the problem
to questions on lower degree curves. This step could be done using
the ``bend and break technique" \cite{Mo}, \cite[II.5]{Ko}, which
gives a deformation $\Gamma '$ of $\Gamma$ whose associated
1-cycle is of the form $\sum a_i\Gamma_i$ with
$\sum a_i\geq 2$ (see \cite[Corollary II.5.6]{Ko}). In order to
fulfil the hypotheses necessary for the application of this
technique, one has to verify that there are enough rational quartics
in $X$ passing through two general points. 
Each component of the Hilbert scheme
of curves in $X$ whose generic point is a smooth rational normal quartic is
of dimension 8. Hence, one can expect that 
in each component, there is a 4-dimensional
family of curves passing through two general points $P_1,P_2\in X$.
To verify that this is really the case, the standard argument with 
incidence varieties can be used (compare to Lemma \ref{lemma-2}).
One has to show that there is a component of dimension 4 in the subset
of the Hilbert scheme of rational quartics passing through two general
points in a {\em particular} cubic threefold. One can easily check
that this is so on a cone over a smooth cubic surface $Y\subset\PP^3$.
\end{remark}

\section{Normal elliptic quintics}

\label{section-quintics}

\subsection{}\label{nota}
Let $H_0$ be the component of $\Hilbp$ whose generic
point corresponds to a smooth normal elliptic quintic in $\PP^4$.
$H$ or $H(X)$ will denote, as before, the union of components of $\Hilb$
having smooth normal elliptic quintics as their generic points
(hence contained in $H_0$).
It is well-known that $H_0$ is unique, of dimension 25, and is generically
parametrized by 24 parameters corresponding to the action of
$\PGL$, and by the modulus of the elliptic curve (see \cite{Hu}, or 
\cite[Proposition 6.1]{HH} for the irreducibility of Hilbert schemes
of curves in a more general context).

\begin{lemma}\label{lemma-1}
Through any smooth elliptic normal quintic in $\PP^4$ passes a $19$-dimensional
linear system of cubic hypersurfaces, and its general member is
nonsingular.
\end{lemma}

\begin{proof}
The dimension calculation follows from (\ref{restriction}) for $k=3$.
By \cite[Theorem IV.1.3]{Hu} or by \cite[Theorem 10]{Mu}, $C$ is a
scheme-theoretic intersection of quadrics passing through $C$.
Hence this is also true for cubics, and by Bertini's Theorem, the
general cubic through $C$ is nonsingular.
\end{proof}

\begin{lemma}\label{lemma-2}
Let $I\subset H_0\times \PP^{34}$ be the incidence variety parametrizing
the pairs $(C,X)$ consisting of a smooth normal 
elliptic quintic $C$ and a cubic hypersurface $X$ containing $C$. 
Then $I$ is irreducible  of dimension 44, and
its projection to $\PP^{34}$ is dominant. 
The fiber of this projection has dimension 10 and is
smooth at any  point $[C]$ representing a normal elliptic quintic $C$
in a nonsingular cubic threefold.
\end{lemma}

\begin{proof}
The irreducibility of $I$ follows from \cite[Theorem 1.6.12]{Sh}. Now, take any
nonsingular cubic $X_0$ containing a smooth normal elliptic quintic $C_0$,
which is possible by Lemma \ref{lemma-1}. Then, by Lemma \ref{h0NCX}, c),
$h^0(\NNN_{C_0/X_0})=10,\; h^1(\NNN_{C_0/X_0})=0$, hence $H(X_0)$ is
smooth of dimension 10 at $[C_0]$ (see \cite{G}). 
\end{proof}

Sofar, we have proved:

\begin{proposition}\label{H-comps}
Let $X$ be a general cubic threefold. Then $H=H(X)$ is the union of
finitely many irreducible components of dimension 10. The generic point of
any of these components is smooth and corresponds to a 
smooth normal elliptic quintic of $\PP^4$ contained in $X$.  
For any such quintic $C$,  the point $[C]\in H$ is smooth.
\end{proposition}

\begin{theorem}\label{H-irred}
Let $X$ be a general cubic threefold, and $H=H(X)$.  Then $H$ is irreducible of
dimension 10.\end{theorem}

\begin{proof}
Choose any component $H'$ of $H$.
By Proposition \ref{H-comps}, there is a smooth normal quintic
$C\subset\PP^4$ corresponding to a point of $H'$. Choose a linear series
$g^1_2$ on $C$, and construct the rational cubic scroll $\Sigma$ as the
union of lines joining pairs of points in $g^1_2$. One can
easily verify that $\Sigma$ is the Hirzebruch surface $\FF_1$
embedded by the linear system of degree 3, and its intersection with
$X$ is $C$ plus a rational normal quartic $\Gamma$ or, possibly, a degenerate
curve of the linear system of rational normal quartics in $\Sigma$.
Such curves $\Gamma$ meet the $(-1)$-curve $L$ of $\Sigma$ in two
points. Notice, that the family of scrolls $\Sigma$ obtained by
this construction is connected of dimension 1, parametrized by  the
points of $\Pic^2(C)\simeq C$.

As a next step, we want to invert this construction and to recover $C$
from $\Gamma$. It is convenient to reduce the problem to the case
when $\Gamma$ is smooth (see Remark \ref{conversely}, giving an idea,
how to avoid this reduction). The singular members of $|\Gamma |$
on $\Sigma$ are of the form $\Gamma'+\sum a_iF_i$, with $1\leq\sum a_i\leq 3$,
where $\Gamma'$ (resp. $F_i$) is a section (resp. a fiber) of the scroll 
$\Sigma$,
and $F_i$ are chords of $C$ lying in $X$. So, if $C$ has only finitely
many chords in $X$, then we can always choose a $g^1_2$ on $C$ in such
a way that it contains no chords, and then $\Gamma$ will be nonsingular.
The finiteness of the number of chords is guaranteed by the following lemma.

\begin{lemma}\label{gen-gen}
A general elliptic quintic $C$ in any component of $H(X)$ for a general
cubic threefold $X$ has only finitely many chords in $X$.
\end{lemma}

Postponing the proof of the lemma, let us show now that
given any normal rational quartic 
$\Gamma\in \PP^4$ contained in $X$,
we can find a cubic scroll $\Sigma$ such that $\Sigma\cap X=\Gamma\cup
C$ for some quintic $C$. It can be constructed geometrically as follows.
Take any chord $L$ of $\Gamma$, meeting $\Gamma$ in two distinct points
$P_1,P_2$.  Choose arbitrarily two points $P_3\in\Gamma$ and $P'_3\in L$,
different from $P_1,P_2$. Connect each point $M\in \Gamma$ by a line
with $M'\in L$, where $M'$ is defined by the equality of the
cross-ratios $[P_1,P_2,P_3,M]=[P_1,P_2,P'_3,M']$. The union of these
lines is a cubic scroll $\Sigma$ having $L$ as $(-1)$-curve. In the case
when $P_1=P_2=P$, $L$ is tangent to $\Gamma$ at $P$, and the scroll
is determined by an isomorphism $\phi :\Gamma\lra L$ such that
$\phi (P)=P$, $d_P\phi =\id_{T_PL}$ (with a natural identification
of $T_PL$ and $T_P\Gamma$).
The residual curve in the intersection with $X$ yields a quintic
$C$.
We see that the family
of the scrolls $\Sigma$ constructed from a fixed quartic
curve is irreducible and has dimension 3. Indeed, 
it is parametrized by a 3-dimensional variety fibered over
$\Sym^2\Gamma$ with fiber over a pair of points $(P_1,P_2)\in\Sym^2\Gamma$
equal to 
\newcommand\Span{\operatorname{Span}\nolimits}
$\Span (P_1,P_2)\smallsetminus\{ P_1,P_2\}$, if $P_1\neq P_2$,
and to the curve of isomorphisms $\Gamma\lra L$ fixing
$P=\Gamma\cap L$ and $T_PL$, if $P_1=P_2=P$.

So, our construction yields a 3-dimensional
irreducible family of scrolls whose residual intersection with $X$
defines a quintic curve. Taking into account Proposition
\ref{H-comps}, we obtain the following statement:

\begin{proposition}\label{fam-quart}
Let $U$ be the open set of $H$ parametrizing general
(in the sense of Lemma \ref{gen-gen}) smooth 
normal elliptic quintics of $\PP^4$ lying on a general
cubic threefold $X$, and $V$
the open subset of  the Hilbert scheme
of quartic curves in $X$ parametrizing all the
curves which are residual intersections with $X$ of the
cubic scrolls through quintics $C$ with $[C]\in U$.
Then the above construction establishes a correspondence
$Z\subset V\times U$ with irreducible fibers of dimension
3 over $V$ and of dimension 1 over $U$. In particular,
$U$ and $V$ have equal numbers of irreducible components.
\end{proposition}

The proof of Theorem \ref{H-irred} is ended by the application
of  Theorem \ref{quart-irred}.
\end{proof}

\subsubsection*{Proof of Lemma \ref{gen-gen}}
Take $C$ smooth and projectively normal in any component
$H_0'\subset H(X_0)$ of a given nonsingular cubic threefold
$X_0$. Let $\pi_P:C\lra C_1$ be a general projection from a point
$P\in\PP^4\setminus X_0$ to $\PP^3\subset\PP^4$. Then $C_1$ is an
elliptic quintic in $\PP^3$. There exists a nonsingular cubic
surface $Y$ containing $C_1$; the proof of this assertion is similar
to that of Lemma \ref{lemma-1}. Let $X_1$ be the cone over $Y$
with vertex $P$. Then $X_1$ is also a cubic threefold containing
$C$. Every chord of $C$ lying in $X_1$ is projected by $\pi$ into
a chord of $C_1$ lying in $Y$. As the number of lines in a nonsingular
cubic surface is equal to 27, the number of chords of $C$ in $X_1$
is finite. Hence it is also finite in a general member $X_\lambda$ of the
pencil of cubics $(X_0,X_1)$. As we can start with $C$ in any
component of $H(X_0)$, we see that a general $C$ in any component
has only finitely many chords in $X_t$ for a general deformation $X_t$
of $X_0$, and we are done.
\hfill\square

\begin{remark}\label{conversely}
The construction of the cubic scroll through a rational quartic
extends also to the case when $\Gamma$
is a singular member of the corresponding linear system on some
scroll $\Sigma\simeq\FF_1$.
Indeed, this linear system is $|L+3F|$, 
where $F$ denotes the fiber of the ruled surface, and the only possible
degenerate members are of one of the following types:
a)  $L+F_1+F_2+F_3$, where $F_i$ are
lines on $X$ meeting $C$ in two points 
(chords), or tangent to $C$; 
b) $C_2+F_1+F_2$, where $C_2$ is a smooth conic;
c) $C_3+F_1$, where $C_3$ is a twisted cubic in some $\PP^3$.
Treat, for example, the case a); the arguments are similar in the
remaining cases. Let $P_i=F_i\cap L$, $i=1,2,3$. 
The three lines, whether they are
distinct or not (in the last case we consider them 
as a non-reduced subscheme of $\PP^4$), do not lie in one hyperplane.
Otherwise, this hyperplane would have the intersection number
at least 6 with $C$, and hence $C$ would lie in the hyperplane,
which contradicts our assumptions. So, through three points
$Q_i\in F_i$ (possibly infinitesimal, in the non-reduced case),
different from $P_i$,
passes a unique 2-plane $\Pi$. Choosing in this plane a conic $B$
through the $Q_i$, we can join the points of $L,B$ with the same
cross-ratios with respect to $P_i=L\cap F_i , Q_i$. This 
construction giving a scroll $\Sigma\supset\Gamma$ works in the
reduced case, and can be extended to the non-reduced one by an
obvious passage to the limit.
\end{remark}

\begin{remark}\label{h1nonnormal}
We saw in the proof of Lemma \ref{lemma-2}, c) that the smooth normal elliptic
quintics are represented by smooth points of the Hilbert scheme of curves.
In fact,
smooth non-normal quintics (that is, lying in a hyperplane), are represented
also by smooth points in $H$, though they may be rather ``singular''
in other situations: they define unstable vector bundles
of rank 2, and their images under the Abel--Jacobi map coincide
with those of pairs of lines $l_1+l_2\; (l_i\subset X)$.

Indeed, $h^0(\NNN_{C/X})=10,h^1(\NNN_{C/X})=0$ by Lemma \ref{h0NCX},
so $[C]_X$ is a smooth point. It belongs to the component $H$ of $\Hilb$,
because the family of elliptic quintics in hyperplane sections of
$X$ is 9-dimensional, so that $C$ has deformations to smooth
quintics not contained in a hyperplane.

\end{remark}

\section{Factorization of the Abel--Jacobi map through moduli of
vector bundles}
\label{factor}

Let $X$ be a general cubic threefold.
Let $M_X(2;0,2)$ be the Gieseker--Maruyama moduli space
\cite{Ma}, \cite{Sim} of semistable (with respect to $\OOO_X(1)$)
rank 2 torsion free sheaves on $X$ with Chern classes $c_1=0$
and $c_2=2[l]$, where $[l]$ is the class of a line modulo
numerical equivalence.
Define the Zariski open subset $M_0$ in it as follows:
\begin{multline}
M_0 =\{[\EEE ]\in M_X(2;0,2) \; |\; (i)\; \EEE\;\mbox{is
{\em stable} and {\em locally free}};\\ (ii) \; h^1(\EEE (-1))=
h^1(\EEE (1))=h^2(\EEE (1))=h^2(\EEE\otimes\EEE )=0
\} .
\end{multline}

\begin{lemma}\label{prop-E}
Let $\EEE$ be a sheaf on $X$ with $[\EEE ]\in M_0$. Then we have:

a) $h^0(\EEE (1))=6,\; h^i(\EEE (1))=0$ for $i>0$.

b) The scheme of zeros $C$ of any section $s\neq 0$ of $\EEE (1)$
is of pure dimension $1$ and is not contained in a hyperplane.

c) $h^i(\EEE (-1))=0\;\; \forall\; i\in\ZZ$.

d) $h^0(\EEE\otimes\EEE )=1$,  $h^1(\EEE\otimes\EEE )=5,
h^2(\EEE\otimes\EEE )=h^3(\EEE\otimes\EEE )=0$.

e) $h^0(\NNN_{C/X})=10,h^1(\NNN_{C/X})=0$.
\end{lemma}

\begin{proof}
This follows immediately from
the definition of $M_0$, the Serre duality, exact sequences
(\ref{EE}), (\ref{timesE1}), and from Riemann--Roch--Hirzebruch
for $\EEE (1),\EEE\otimes\EEE$.
\end{proof}

Define now the locally closed subscheme $\HHH_0\subset \Hilb$:

\begin{equation}\label{defH0}
\HHH_0=\left\{\begin{minipage}{8.8 truecm}
$[C]\in\Hilb\; |\;$ (i) $C$ is a locally complete intersection of pure
dimension 1,
(ii) $h^1(\III_C)=h^1(\III_C(1))=0$ (hence $h^0(\OOO_C)=1$),
(iii) $h^1(\III_C(2))=h^2(\III_C(2))=0$ (hence $h^0(\III_C(2))=5$),
(iv) $\omega_C\simeq \OOO_C$
\end{minipage}\right\}
\end{equation}

\begin{lemma}[Relative Serre's construction]\label{lem61}
There is a well-defined \linebreak morphism
$\phi :\HHH_0\lra M_0,\; [C]\mapsto [\EEE ]$,
where $\EEE$ is defined by the exact triple (\ref{serre}),
determined by a non-zero extension class in $H^0(\omega_C)
=\Ext^1(\III_C(2),\OOO_X)$.
\end{lemma}

\begin{proof}
The proof consists in an obvious relativization of Serre's construction
of Sect. \ref{sect.5} (see Corollary \ref{unique-ext}).

\end{proof}

\begin{lemma}\label{lem62}
$\HHH_0$ is isomorphic to
the projectivization of a rank 6 vector bundle
locally in the \'etale topology over $M_0$.
\end{lemma}

\begin{proof}
According to Simpson \cite[Theorem 1.21]{Sim}, each point of $M_0$
has an \'etale neighborhood $T\lra U\subset M_0$, such that there exists
a Poincar\'e vector bundle $\EEEE$ on $X\times T$. Let
$\GGGG =\pr_{2*}(\EEEE\otimes \OOO_X(1))$, a locally free
sheaf of rank 6 on $T$.
Let $\CCC$ be the universal family of curves over
$\HHH_0$, and $\CCC_U\lra\HHH_0^U$ its restriction over $U$. 

Given a commutative diagram
$$
\begin{CD}
Y @>{\alpha}>> \HHH_0^U \\
@V{\beta}VV @VV{\phi}V \\
T @>{\lambda}>> U
\end{CD}
$$
we will show that there exists a unique morphism $\gamma :Y\lra\PP (\GGGG )$
such that $q\circ\gamma =\alpha ,\; p\circ\gamma =\beta$, 
where $p:\PP (\GGGG )\lra T$ is the structure morphism and
$q:\PP (\GGGG )\lra \HHH_0$ is the natural map sending the proportionality
class of a section $s_t$ of $\EEE_t(1)$ on $X\times \{ t\}$ to
its scheme of zeros $Z_X(s_t)$. Then, by the universal property of the
Cartesian product, $\PP (\GGGG )=\HHH_0^U\times_U T$.
 
Lift $\CCC_U$ to the family $p_2: \CCC_Y\lra \HHH_0^Y=Y\times_U\HHH_0$.
By conditions (ii), (iv) of 
(\ref{defH0}) and Base Change, $p_{2*}\omega_{\CCC_{Y}/Y}$ is
invertible, so that there exists an open covering
$Y=\cup_{i\in I}Y_i$  and non-vanishing sections $\xi_j
\in\Gamma (Y_j, p_{2*}\omega_{\CCC_Y/Y})$. They determine extensions
$$
0\lra\OOO_{X\times Y_j}
\stackrel{\mu_j}{\lra}\tilde{\EEEE}_{j}(1)\lra\III_{\CCC_j}(2)\lra 0 ,
$$
where $\CCC_j=\CCC_Y|_{Y_j}$, such that
$[\tilde{\EEEE}_{j}|_{X\times\{ y\} }]=\lambda\circ\beta (y)\;
\forall\; y\in Y_j$.

Hence, denoting $\beta_j=\beta|_{Y_j}$, we get by the universal property
of the Poincar\'e bundle the isomorphisms
$\tilde{\EEEE}_{j}\simeq (1\times\beta_j)^*\EEEE\otimes \MMM_j$ for
some invertible sheafs $\MMM_j$  on $Y_j$. Applying
$p_{2*}$ to $\mu_j$, we obtain the monomorphism
$p_{2*}\mu_j:\OOO_{Y_j}\hookrightarrow \beta_j^*\GGGG\otimes \MMM_j$,
or, by duality, an epimorphism $\epsilon_j:\beta_j^*\GGGG^*
\lra\MMM_j$. By the universal property of the functor
{\bf Proj} \cite[Exercise II.7.8]{Ha}, there exists a unique morphism
$\gamma_j:Y_j\lra\PP(\GGGG )$ such that $\epsilon_j=\gamma_j^*\epsi$,
where $\epsi :p^*\GGGG^*\lra\OOO_{\PP(\GGGG )/T} (1)$ is the canonical
epimorphism.
These morphisms agree on the intersections $Y_i\cap Y_k$ and 
thus define a morphism $\g :Y\lra \PP(\GGGG )$. By construction, 
$q\circ\gamma =\alpha ,\; p\circ\gamma =\beta$.
\end{proof}

\begin{corollary}\label{sharp}
The morphism $\phi :\HHH_0\lra M_0$ defined in Lemma \ref{lem61}
is smooth, projective, and all its fibers are $5$-dimensional
projective spaces.
\end{corollary}

\begin{proof}
According to Lemma \ref{lem62}, $\phi :\HHH_0\lra M_0$ is identified,
locally in the \'etale topology, with the natural projection
$\PP^5\times U\lra U$.
\end{proof}

\begin{corollary}\label{flat}
$\HHH_0,M_0$ are smooth of dimensions $10$, resp. $5$; moreover, $\Hilb$,
$M_X(2;0,2)$ are smooth at the points of $\HHH_0$, resp. $M_0$.
\end{corollary}

\begin{proof}
This follows from Lemma \ref{prop-E} d), e).
\end{proof}

We do not treat the question of irreducibility of $\HHH_0$ and $M_0$.
Let $\HHH_0'$ be the irreducible component of $\HHH_0$ containing the open
subset
$$
\HHH^*=\{ [C]\in \HHH_0\; |\; C\;\mbox{is a smooth projectively normal
elliptic quintic}\} .
$$
$\HHH^*$ is irreducible by Theorem  \ref{H-irred},
and we proved in Section \ref{sect.5} that $\HHH^*\subset \HHH_0$.
Define also $M=\phi (\HHH^*),\HHH =\phi^{-1}(M)\cap\HHH_0'$. In the sequel,
we will denote the restriction $\phi |_{\HHH}$ by the same symbol $\phi$.
By Corollaries \ref{sharp},
\ref{flat}, $M,\HHH$ are open and smooth.

\begin{theorem}\label{main}
Let $\HHH^*\subset\HHH\subset\Hilb , M\subset M_X(2;0,2)$, $\phi :\HHH
\lra M$ be defined as above. The following assertions are true:

(i) For every choice of a reference point $[C_0]\in\HHH$, the Abel--Jacobi map
$[C]\mapsto [C-C_0]$ defines a morphism $\Phi :\HHH\lra J_1(X)$.
Fix some $[C_0]$ and the corresponding morphism.

(ii) $\Phi$ is smooth, and every fiber of $\Phi$ is a disjoint union
of $5$-dimensional projective spaces.

(iii) There exists a quasi-finite \'etale morphisme $\Psi :M\lra J_1(X)$
such that $\Phi=\Psi\circ\phi$.
\end{theorem}

\begin{proof}
(i) $\HHH$ is smooth and is a base of a flat family of curves,
hence, by \ref{AJ_B}, the Abel--Jacobi map $\Phi :\HHH\lra J_1(X)$
is well-defined as an analytic map. It is in fact algebraic. Indeed,
$\HHH$ is a subvariety of the Chow variety of $1$-cycles in $X$
(see \cite[Theorem I.6.3]{Ko}), and we can resolve the singularities
of the closure of $\HHH$ in the Chow variety. Then $\Phi$ extends to this
resolution as an analytic map. It will be projective by GAGA principle,
as the Chow variety is projective.

(ii) First, we will verify that $\Phi$ is smooth at any point of $\HHH^*$.
This follows from the computation of the differential of $\Phi$ by the
technique
of TBS (see \ref{TBS}). Apply it to $X\subset W=\PP^4$. The map $R$ of
(\ref{CDWelters}) is an isomorphism $H^0(X,\OOO_X(1))
\stackrel{\textstyle\sim}{\lra}H^1(X,\Omega_X^2)$ by \cite[(12.7)]{CG}.
We have to verify that the map $\psi_Z^*$ of (\ref{CDWelters}) is injective
for $Z=C$ with $[C]\in \HHH^*$. It suffices to show the injectivity of
$r_C,\beta_C$. The kernel of $r_C$ lies in
$H^0(X,\NNN_{X/\PP^4}\otimes\omega_X\otimes\III_{C/X})=
H^0(\III_C(1))=0$ by Lemma \ref{dim=6} b). By (\ref{CDsheaves}),
$\beta_C$ is a part of the exact sequence
\begin{multline*}
\ldots\lra H^0(\Omega_X^3\otimes\NNN_{C/\PP^4})\lra
H^0(\Omega_X^3\otimes\NNN_{X\PP^4}|_C)\stackrel{\beta_C}{\lra} \\
H^1(\Omega_X^3\otimes\NNN_{C/X})=H^0(\NNN_{C/X})^* .
\end{multline*}
By \cite[V.2.1]{Hu}, $H^0(\Omega_X^3\otimes\NNN_{C/\PP^4})=0$.
Hence $\beta_C:H^0(\OOO_C(1))\lra $\linebreak $H^0(\NNN_{C/X})^*$ is injective.
Hence $\Phi$ is of maximal rank, equal to 5, at $[C]$.

Thus, every $z\in \HHH^*$ is contained in a 5-dimensional component
$\HHH_z$ of the fiber $\Phi^{-1}\Phi (z)$, nonsingular at $z$.
We know also that every projective space contained in $\HHH$ is contracted 
by $\Phi$, because the Abel--Jacobi map contracts the rational 
equivalence. Hence $\HHH_z=\phi^{-1}\phi (z)$, the fiber of~$\phi$.

More generally, every fiber $\Phi^{-1}(w)$ is the union of projective
spaces $\phi^{-1}\phi (z)\simeq\PP^5$ over all $z\in\Phi^{-1}(w)$,
and for $z\neq z'\in\Phi^{-1}(w)$ we have: either $\phi^{-1}\phi (z)=
\phi^{-1}\phi (z')$, or $\phi^{-1}\phi (z)\cap
\phi^{-1}\phi (z')=\varnothing$.
By construction, every fiber $\PP^5$ of $\phi$ contains a point of
$\HHH^*$, hence $\Phi^{-1}(w)$, if non-empty, is a disjoint finite union of
irreducible components of the form $\phi^{-1}\phi (z)\simeq\PP^5$
with $z\in\Phi^{-1}(w)$, and there are no multiple components among them,
since the fiber should be scheme-theoretically smooth at points
of $\HHH^*$. Moreover, the fibers of $\Phi$ are locally complete intersections
in $\HHH$, because they are of pure codimension
5 and $J_1(X)$ is smooth of dimension 5, so, they have
no embedded points. Hence
all the components of the fibers of $\Phi$ are smooth
subvarieties of $\HHH$, isomorphic to $\PP^5$, with their
reduced structure.

(iii) By (ii), the map $\Psi$ is well-defined set-theoretically:
we saw that $M$ parametrizes the irreducible components (isomorphic
to $\PP^5$) of the fibers of $\Phi$. By Corollary \ref{sharp},
$\phi_*\OOO_\HHH =\OOO_M$, and $\phi$ is open. For every open $U\subset
J_1(X)$, the inverse image $\Psi^{-1}(U)$ under our set-theoretical
map is open, because $\Psi^{-1}(U)=\phi(\Phi^{-1}(U))$, and $\Psi^*(g)$
is a regular function on $\Psi^{-1}(U)$ for every $g\in\OOO_{J_1(X)}(U)$.
Thus $\Psi$ is Zariski continuous and lifts regular functions on
a Zariski open subset $U\subset J_1(X)$ to those on $\Psi^{-1}(U)$.
Hence $\Psi$ is a morphism.
We have $d\Phi =d\Psi\circ d\phi$ and $d\Phi$ is surjective. Hence
$d\Psi$ is surjective, hence it is an isomorphism and $\Psi$ is \'etale.
\end{proof} 

Leaving aside the hard question on the complete description of the boundary
of $\HHH$ in $\Hilb$ and of its image in $J_1(X)$, we can easily determine
the image of one of the components of the boundary, namely, the one
whose general point represents a projectively non-normal elliptic quintic.
Let $H_\ns$ be the nonsingular locus of $H$,  $B^*\subset H$
the subset of classes of smooth elliptic quintics contained in
nonsingular hyperplane sections of $X$, and $B$ the closure of $B^*$ in $H$.
By Remark \ref{h1nonnormal}, $B^*\subset H_\ns$ and $B$ are of pure dimension
9. According to \ref{AJ_B}, the Abel--Jacobi map $\Phi$ of Theorem \ref{AJ_B}
extends to a morphism, well defined on $H_\ns$, hence on $B^*$ and on
the desingularization of $B$.

\begin{proposition}\label{aftermain}
Let $\tilde{\Phi}:H_\ns \lra J_1(X)$ be the Abel--Jacobi map defined
above. The following properties are verified:

(i) $B$ is an irreducible divisor in $H$.

(ii) The image $\tilde{\Phi}(B^*)$ of $B^*$ in $J_1(X)$
is dense in the translate $\Delta_a$ by an element $a\in J_1(X)$ 
of the divisor $\Delta =\tilde{\psi}
(F\times F)$ defined in \ref{phi}.

(iii) Let $\eta =\Phi |_{B^*}:B^*\lra\Delta_a$ be the restriction
of $\tilde{\Phi}$. Then the general fiber of $\eta$ is an open subset
of the 5-dimensional projective space $\PP^5$, realized as the
subset of nonsingular curves in a linear system on a cubic surface.

(iv) $\tilde{\Phi}$ is ramified along $B$. More exactly,
$\corank \: d\tilde{\Phi}\: =\: 1$ at every point of $B^*$.
\end{proposition}

\begin{proof}
By Proposition \ref{space-C}, a curve $C$ with $[C]\in B^*$
determines a unique unordered pair of disjoint lines $l_1,l_2$
on the nonsingular cubic surface
$S=X\cap <C>$ in such a way that $C\in |h+l_1+l_2|\simeq\PP^5$, where
$<C>$ is the linear span of $C$ in $\PP^4$, and $h$ is the hyperplane
section of $S$. The pairs of disjoint lines are parametrized by an irreducible
4-fold $\Sym^2F\setminus\mbox{(incidence divisor)}$, and the pairs whose
span defines a nonsingular hyperplane section form an open subset
$U$ in it. Hence $B^*$ is irreducible,
and its image in $J_1(X)$ coincides with $\tilde{\psi}(U)$
modulo a translation depending
on the choice of reference points for $\tilde{\psi},\Phi$. 
By Beauville's remark in \ref{beau}, the degree of $\tilde{\psi}$
over its image is 2, so there is a unique unordered pair $l_1,l_2$ over
the generic point of $\Delta_a$. This proves 
(i)--(iii). The assertion (iv) follows from the technique of TBS: as in the
proof of Theorem \ref{main}, (ii), for $[C]\in B^*$, we find that 
$\beta_C$ is injective, but $\ker\: r_C=H^0(\III_C(1))$ is of dimension
1, which implies the result.
\end{proof}

\end{document}